# A pattern mixture model for a paired $2 \times 2$ crossover design

**Laura J. Simon**[1] **and Vernon M. Chinchilli**[*,2]

*The Pennsylvania State University*

**Abstract:** When conducting a paired $2 \times 2$ crossover design, each subject is paired with another subject with similar characteristics. The pair is then randomized to the same sequence of two treatments. That is, the two subjects receive the first experimental treatment, and then they cross over and receive the other experimental treatment(s). The paired $2 \times 2$ crossover design that was used in the Beta Adrenergic Response by GEnotype (BARGE) Study conducted by the National Heart, Lung and Blood Institute's Asthma Clinical Research Network (ACRN) has been described elsewhere. When the data arising from such a design are balanced and complete – or if at least any missingness that occurs is at random – general linear mixed-effects model methods can be used to analyze the data. In this paper, we present a method based on a pattern-mixture model for analyzing the data arising from a paired $2 \times 2$ crossover design when some of the data are missing in a non-ignorable fashion. Because of its inherent scientific interest, we focus our particular attention on the estimation of the treatment-by-type of subject interaction term. Finally, we illustrate the pattern-mixture model methods described in this paper on the data arising from the BARGE study.

## 1. Introduction

Two important design principles are occasionally used in clinical trials: 1) A subject is "matched" or "paired" with another subject with similar characteristics to reduce the chance that other variables obscure the primary comparison of interest. 2) A subject serves as his or her own control by "crossing over" from one treatment to another during the course of an experiment. Experiments employing the first design principle are called matched pair designs (Cochran [1]), while those employing the second are called crossover designs (Jones and Kenward [2], Ratkowsky et al. [3] and Senn [4]).

There are situations in which it may be beneficial to use the two design principles simultaneously. That is, it may be advantageous to conduct a "paired crossover design." For such a design, each subject and his/her paired counterpart are randomized to the same treatment sequence. That is, they receive one experimental treatment, and then cross over and receive the other experimental treatment(s) at the same time. A paired crossover design was recently used in three clinical tri-als conducted by the National Heart, Lung and Blood Institute's Asthma Clinical Research Network (ACRN) (Kephart [5]; www.acrn.org).

*Supported in part by Grants U10 HL51845 and U10 HL074231 from the National Heart, Lung, and Blood Institute

[1]Department of Statistics, The Pennsylvania State University, 328 Thomas Building, University Park, PA 16802, USA, e-mail: `lsimon@stat.psu.edu`
[2]Department of Public Health Sciences, The Pennsylvania State University, A210, 600 Centerview Drive, Suite 2200, Hershey, PA 17033, USA, e-mail: `vchinchi@psu.edu`

In particular, the Beta Adrenergic Response by GEnotype (BARGE) Study (Israel et al. [6]) used a paired $2 \times 2$ crossover design to compare the regular use of inhaled albuterol (A) to placebo (P) in patients with the Arg/Arg genotype and the Gly/Gly genotype at the $16^{th}$ position of the beta-agonist receptor gene. It was hypothesized that the scheduled daily use of inhaled albuterol, the most common treatment for patients with mild to moderate asthma, actually had a detrimental effect on the lung function of patients with the Arg/Arg genotype (R) but not patients with the Gly/Gly genotype (G). The primary research question therefore concerned whether or not the treatment effects differed for the two genotypes. That is, the primary hypothesis concerned inference about whether the interaction parameter:

$$\gamma = (\mu_{RA} - \mu_{RP}) - (\mu_{GA} - \mu_{GP})$$

is 0, where $\mu_{kl}$ is the population mean lung function of genotype $k$ patients on treatment $l$.

To achieve an efficient comparison of the treatments within each genotype, a $2 \times 2$ crossover design was used within each genotype $k$. Because a subject's genotype is a pre-determined characteristic, subjects could not be randomly allocated to genotypic group. Therefore, to minimize confounding caused by potential differences in the baseline lung function of the two genotypic groups, each subject of genotype R was matched to a subject of genotype G with similar baseline lung function. The matched subjects were randomly assigned to the same sequence of the crossover design, and an eight-week washout period was placed between the two treatment periods (Table 1).

At the conclusion of the BARGE Study, each pair of subjects $j$ of sequence $s$ ideally yielded the quadrivariate response $\mathbf{Y}_{sj} = (Y_{sjRA}, Y_{sjRP}, Y_{sjGA}, Y_{sjGP})$, a vector containing the subjects' changes in lung function. Unfortunately, as should be expected when conducting any clinical trial, measurements in the BARGE Study were not always collected as planned. That is, some subjects missed one or more planned visits or completely dropped out of the study. In this paper, we present a method based on a pattern-mixture model for analyzing the data arising from a paired $2 \times 2$ crossover design when some of the data are missing.

In general, if $\mathbf{y}_i = (y_{i1}, \ldots, y_{ip})'$ denotes the $p \times 1$ vector of intended repeated measurements on subject $i$ and $\mathbf{m}_i = (m_{i1}, \ldots, m_{ip})'$ represents a $p \times 1$ vector of random indicator variables denoting whether or not response $y_{ij}$ is observed, then pattern-mixture models factor the joint distribution of $\mathbf{y}_i$ and $\mathbf{m}_i$ as:

$$\pi\left[\mathbf{y}_i, \mathbf{m}_i | \mathbf{X}_i\right] = \pi\left[\mathbf{y}_i | \mathbf{X}_i, \mathbf{m}_i\right] \times \pi\left[\mathbf{m}_i | \mathbf{X}_i\right],$$

where $\mathbf{X}_i$ denotes fixed covariates or design matrices. That is, in short, the data are stratified by their patterns of missingness, and then a separate model is specified for each missing data pattern. The distribution $\pi[\mathbf{y}_i|\mathbf{X}_i, \mathbf{m}_i]$ models the within-subject regressions for each missing data pattern, and $\pi[\mathbf{m}_i|\mathbf{X}_i]$ models the marginal proportions of each missing data pattern as functions of between-subject covariates.

TABLE 1
*Paired crossover design for the BARGE Study. Pairs of subjects were randomized to either the AP sequence (first row) or the PA sequence (second row). All subjects had a washout period between the two treatment periods*

|            | Subject 1 (R) |    | Washout | Subject 2 (G) |    |
|------------|---------------|----|---------|---------------|----|
| **Period** | 1             | 2  | ——      | 1             | 2  |
| Sequence 1 | RA            | RP | ——      | GA            | GP |
| Sequence 2 | RP            | RA | ——      | GP            | GA |

It should be noted that there does exist a more general random-coefficient pattern-mixture model, e.g., see Little [7]. We do not pursue the more general model here, however, because our crossover model for the BARGE Study assumes that only one measurement is made in each period.

The type of "missingness" that exists really should dictate our final analysis. Rubin [8] and Little and Rubin [9] describe three types of missingness for which we need not apply a pattern-mixture model: i) if the data are missing completely at random (the probability of response is independent of both the observed and unobserved data); ii) if covariate-dependent dropout exists (the missingness depends on fixed covariates in the model); and iii) if the data are missing at random (the probability of response depends on the observed data but not the missing data). In any of these three cases, we simply can apply a general linear model with correlated errors to the data arising from a paired crossover design, as was performed by Simon and Chinchilli [10]. The resulting likelihood-based estimation and inference procedures are asymptotically unbiased.

If the missing values satisfy none of the three types mentioned above, then Little and Rubin [9] label the situation "non-ignorable" (the probability of response depends on the unobserved data), and the analysis of the available data require special methods. In this situation, when analyzing the data arising from a paired crossover design, we propose applying the pattern-mixture model methods that we now describe.

## 2. Methods

In defining a pattern-mixture model for a paired $2 \times 2$ crossover experiment, we assume that each subject is scheduled to contribute one response per period. The response may truly be a single post-treatment, outcome measure, such as forced expiratory volume in one second ("$FEV_1$") in an asthma trial. Alternatively, and perhaps more commonly, the response could be a summary of multiple, repeated outcome measures, such as the change in response from pre- to post-treatment or the area under a dose-response curve ("$AUC$") in a bioequivalence trial. We also assume that the subjects do not contribute any additional covariates. Thus, in general, if we refer to the two types of subjects who are matched as 1 and 2 and we label the treatments as A and B, then we expect each pair of subjects $j$ of sequence $s$ to yield only the quadrivariate response $\mathbf{Y}_{sj} = (Y_{sj1A}, Y_{sj1B}, Y_{sj2A}, Y_{sj2B})$. To define a pattern-mixture model for a paired $2 \times 2$ crossover design that accommodates non-ignorable missing data, we must first consider the possible patterns of missingness.

### 2.1. Patterns of missingness

When modeling missing data in a clinical trial, it is common to assume that the missingness happens monotonically, *i.e.*, after the first missing data point for a subject, all of the subsequent anticipated data points for that subject are also missing. For example, the data vector for a subject $i$ might look something like $\mathbf{y}_i = (y_{i1}, y_{i2}, y_{i3}, \cdot, \cdot, \cdot)'$ assuming the subject misses all measurement occasions after the third occasion. The monotonic missingness assumption greatly simplifies the model development and subsequent analysis.

Unfortunately, for a paired crossover design, the monotonicity assumption is not realistic. The data vector $\mathbf{y}_{sj}$ for the paired $2 \times 2$ crossover design contains data on the two subjects within a pair. Therefore, data vectors with non-monotonic missingness that look like $\mathbf{y}_{sj} = (y_{sj1A}, \cdot, y_{sj2A}, \cdot)'$ are quite possible in a paired crossover design. One could still hope to simplify the problem by assuming that at least the missingness with respect to each subject is monotonic. However, the BARGE Study data set, to which we later apply the methods described herein, contains non-monotonic missingness within subjects. Therefore, we do not make any simplifying assumptions, but instead address the problem in its full generality.

Another complication of handling missing data in a crossover design (paired or otherwise) is that, even when the same pattern of missingness occurs, different information is gleaned about the treatments from the subjects in different sequences. For example, for our paired $2 \times 2$ crossover design, data vectors for the first sequence, $\mathbf{y}_{1j}$, which are missing the last data point:

$$\mathbf{y}_{1j} = (y_{1j1A}, y_{1j1B}, y_{1j2A}, \cdot)'$$

provide no information about the type 2-treatment B combination, while data vectors for the second sequence, $\mathbf{y}_{2j}$, which are missing the last data point:

$$\mathbf{y}_{2j} = (y_{2j1B}, y_{2j1A}, y_{2j2B}, \cdot)'$$

provide no information about the type 2-treatment A combination. Therefore, pattern-mixture models for any crossover design must accommodate patterns of missingness for each of the sequences.

Let $P_{ps}$ denote the pattern of missingness for pattern $p = 0, 1, \ldots, 14$ for sequence $s = 1, 2$. Table 2 summarizes, for the paired $2 \times 2$ crossover design, the eight potential patterns of missingness that are monotonic within subject, while Table 3 summarizes the remaining seven patterns of missingness that are not monotonic within subject.

TABLE 2
*Eight patterns of missingness in paired 2-by-2 crossover design that are monotonic within subject. X denotes observed and ? denotes missing*

|         |          | Subject 1 (R) |       | Subject 2 (G) |       |
|---------|----------|---------------|-------|---------------|-------|
| Pattern | Sequence | Per 1         | Per 2 | Per 1         | Per 2 |
| $P_{01}$ | 1 (AB)  | X | X | X | X |
| $P_{02}$ | 2 (BA)  | X | X | X | X |
| $P_{11}$ | 1 (AB)  | X | X | X | ? |
| $P_{12}$ | 2 (BA)  | X | X | X | ? |
|         |          |   |   |   |   |
| $P_{21}$ | 1 (AB)  | X | ? | X | X |
| $P_{22}$ | 2 (BA)  | X | ? | X | X |
| $P_{31}$ | 1 (AB)  | X | ? | X | ? |
| $P_{32}$ | 2 (BA)  | X | ? | X | ? |
|         |          |   |   |   |   |
| $P_{41}$ | 1 (AB)  | X | X | ? | ? |
| $P_{42}$ | 2 (BA)  | X | X | ? | ? |
| $P_{51}$ | 1 (AB)  | ? | ? | X | X |
| $P_{52}$ | 2 (BA)  | ? | ? | X | X |
|         |          |   |   |   |   |
| $P_{61}$ | 1 (AB)  | X | ? | ? | ? |
| $P_{62}$ | 2 (BA)  | X | ? | ? | ? |
| $P_{71}$ | 1 (AB)  | ? | ? | X | ? |
| $P_{72}$ | 2 (BA)  | ? | ? | X | ? |

TABLE 3
*Seven patterns of missingness in paired 2-by-2 crossover design that are non-monotonic within subject. X denotes observed and ? denotes missing*

| Pattern | Sequence | Subject 1 (R) Per 1 | Per 2 | Subject 2 (G) Per 1 | Per 2 |
|---|---|---|---|---|---|
| $P_{81}$ | 1 (AB) | ? | ? | ? | X |
| $P_{82}$ | 2 (BA) | ? | ? | ? | X |
| $P_{91}$ | 1 (AB) | ? | X | ? | ? |
| $P_{92}$ | 2 (BA) | ? | X | ? | ? |
| $P_{10,1}$ | 1 (AB) | ? | X | X | X |
| $P_{10,2}$ | 2 (BA) | ? | X | X | X |
| $P_{11,1}$ | 1 (AB) | X | X | ? | X |
| $P_{11,2}$ | 2 (BA) | X | X | ? | X |
| $P_{12,1}$ | 1 (AB) | ? | X | ? | X |
| $P_{12,2}$ | 2 (BA) | ? | X | ? | X |
| $P_{13,1}$ | 1 (AB) | ? | X | X | ? |
| $P_{13,2}$ | 2 (BA) | ? | X | X | ? |
| $P_{14,1}$ | 1 (AB) | X | ? | ? | X |
| $P_{14,2}$ | 2 (BA) | X | ? | ? | X |

Complicating issues arise when modeling the data from a $2 \times 2$ crossover experiment even when the data are complete. Therefore, before defining a pattern-mixture model that accommodates non-ignorable missing data, we first define a statistical model assuming the data are complete.

## 2.2. A paired $2 \times 2$ crossover model for complete data

In defining a statistical model for the paired $2 \times 2$ crossover design, one of the complicating issues is the handling of "treatment carryover." In an AB sequence, a subject receives treatment A in the first period and treatment B in the second period. In this case, the "A carryover" ($\lambda_A$) is the component of the response in the second period due to the lasting effect of treatment A from the first period (Jones and Kenward [2], Ratkowsky et al. [3] and Senn [4]). Likewise, in a BA sequence, the "B carryover" ($\lambda_B$) is the component of the response in the second period due to the lasting effect of treatment B from the first period. In a $2 \times 2$ crossover, if the magnitude of the A carryover differs from that of the B carryover ($\lambda_A \neq \lambda_B$), it is not possible to estimate the primary quantity of interest – the true treatment effect $\mu_A - \mu_B$.

There are two situations in which it is possible to estimate the true treatment effect: 1) if the carryover effects are the same for each treatment ($\lambda_A = \lambda_B$); or 2) if the time periods between the treatment periods, i.e., the washout periods, are designed to be lengthy enough to render both carryover effects negligible ($\lambda_A = \lambda_B = 0$). Most $2 \times 2$ crossover experiments – the BARGE Study included – institute a substantial washout period between the treatment periods with the intent of eliminating any carryover effects. We therefore proceed in formulating a model assuming the second situation holds. In general, if unequal carryover effects are suspected to exist, then a design more complex than the $2 \times 2$ crossover is needed.

A general linear model with correlated errors – which is a special case of the general linear mixed-effects model described by Laird and Ware [11] – for the paired $2 \times 2$ crossover design with null carryover effects can be applied and is stated as:

$$\mathbf{Y}_{sj} = \mathbf{X}_{sj}\beta + \varepsilon_{sj},$$

where:

- $\mathbf{Y}_{sj} = (Y_{sj1A}, Y_{sj1B}, Y_{sj2A}, Y_{sj2B})'$ is the quadrivariate response for pair $j$ $(1, 2, \ldots, n_s)$ of sequence $s$ $(1 = AB, 2 = BA)$.
- $\mathbf{X}_{sj}$ is a $4 \times 8$ fixed-effects design matrix (in a complete data situation, $\mathbf{X}_{1j}$ and $\mathbf{X}_{2j}$ do not depend on $j$, i.e., $\mathbf{X}_{sj} = \mathbf{X}_s$ for all $s = 1, 2$ and $j = 1, 2, \ldots, n_s$).
- $\beta = (\mu_{1A}, \mu_{1B}, \mu_{2A}, \mu_{2B}, \rho_1, \rho_2, \nu_1, \nu_2)'$ is an $8 \times 1$ fixed-effects parameter vector containing four type-by-treatment means, two type-by-period effects, and two type-by-sequence effects.
- $\varepsilon_{sj} = (\varepsilon_{sj1A}, \varepsilon_{sj1B}, \varepsilon_{sj2A}, \varepsilon_{sj2B})'$ is a random error term.

The mean of the responses derives directly from the $X_{sj}\beta$ portion of the model, while the variance of the responses derives directly from the $\varepsilon_{sj}$ portion of the model.

### 2.2.1. Mean of the responses

Our proposed model for the mean of the responses arising from a paired $2 \times 2$ crossover design is a straightforward extension of the model commonly assumed for the mean of the responses arising from a basic $2 \times 2$ crossover design. That is, we assume that the mean of the response in any cell of the paired $2 \times 2$ crossover design is a function of the subject type and treatment ($\mu_{kl}$), as well as the period ($\rho_k$) and sequence ($\nu_l$) in which the treatments are taken. More specifically, we propose parameterizing the mean of the responses as:

| Type | Sequence | Period 1 | Period 2 |
|------|----------|----------|----------|
| 1 | AB | $\mu_{1A} + \rho_1 + \nu_1$ | $\mu_{1B} - \rho_1 + \nu_1$ |
| 1 | BA | $\mu_{1B} + \rho_1 - \nu_1$ | $\mu_{1A} - \rho_1 - \nu_1$ |
| 2 | AB | $\mu_{2A} + \rho_2 + \nu_2$ | $\mu_{2B} - \rho_2 + \nu_2$ |
| 2 | BA | $\mu_{2B} + \rho_2 - \nu_2$ | $\mu_{2A} - \rho_2 - \nu_2$ |

### 2.2.2. Variance of the responses

Traditionally, the variance matrix for a crossover design is assumed to be compound symmetric. That is, the variance of the responses arising from one treatment is assumed to equal the variances of the responses arising from other treatments. And, the covariances between the responses of the subjects when on different treatments are also assumed to be equal. Like the approach of others (Ekbohm and Melander [12], Sheiner [13], Chinchilli and Esinhart [14] and Putt and Chinchilli [15]), we instead model the maximum number of variance components permitted for the paired $2 \times 2$ crossover design. That is, the variance, $\boldsymbol{\Sigma}$, of the response $\mathbf{Y}_{sj}$, $s = 1, 2$:

$$\begin{bmatrix} \sigma_{1A,1A} & \sigma_{1A,1B} & \sigma_{1A,2A} & \sigma_{1A,2B} \\ \sigma_{1A,1B} & \sigma_{1B,1B} & \sigma_{1B,2A} & \sigma_{1B,2B} \\ \sigma_{1A,2A} & \sigma_{1B,2A} & \sigma_{2A,2A} & \sigma_{2A,2B} \\ \sigma_{1A,2B} & \sigma_{1B,2B} & \sigma_{2A,2B} & \sigma_{2B,2B} \end{bmatrix}$$

is much more flexible than the compound symmetric structure traditionally assumed.

*2.2.3. The data*

As the paired $2 \times 2$ crossover model for the complete data situation suggests, we assume that the data are balanced (the occasions of measurement are the same for all subjects) and complete (measurements are available at each planned occasion for each subject). However, the model is appropriate even when the data are not complete but the missingness occurs at random (the probability of response depends on the observed data but not the missing data). We now define a pattern-mixture model for when the missingness is non-ignorable.

**2.3. A pattern-mixture model**

Attempting to define a model for each of the fifteen missing data patterns yields identifiability problems. Some of the patterns will naturally be more populated, while others will be sparsely populated at best. For the sake of illustration, Table 4 summarizes the number of subject pairs in the BARGE Study falling into each of the fifteen missing data patterns. Table 4 indicates that there were no non-monotonic patterns observed in the BARGE Study, which is not suprising because of the longitudinal nature of the trial design.

In general, because of the inherent identifiability problems, it is necessary to collapse the fifteen patterns into coarser groupings, so that information about the effects parameters can be "borrowed" across the patterns. When considering potential groupings, one should take into account the area of scientific research, as well as the various characteristics of the subjects and/or pairs that might yield group differences. In creating groupings for the BARGE Study data, we propose the formation of the following three groups:

1. A group containing the completers who have a pair match – patterns 0, 10, 11, and 12. This group is denoted group "**C**" for "completers."
2. A group in which all patterns involve missing a second period, regardless of whether a pair match exists – patterns 1, 2, 3, 6, 7, 13 and 14. This group is denoted group "**D**" for missingness due to "dropout."
3. A group of the patterns in which a subject is a completer but does not have a pair match – patterns 4, 5, 8 and 9. This group is denoted group "**P**" for missingness due to a missing "pair."

Collapsing the counts in Table 4 for the BARGE Study according to these criteria yields that $n_C = 29$, $n_D = 6$, and $n_P = 5$.

It should be emphasized that other groupings are possible besides the one that we propose. Some of the other groupings, however, require one to assume that no period and no sequence effects exist. That is an assumption that we were unwilling to make for the BARGE Study. It should be noted, though, that if there were no period nor sequence effects to worry about, the methods described by Little (1995) could be applied instead.

TABLE 4
*Frequencies of patterns observed in the BARGE Study*

| Pattern | Count | Pattern | Count | Pattern | Count |
|---|---|---|---|---|---|
| 0 | 29 | 5 | 2 | 10 | 0 |
| 1 | 1 | 6 | 1 | 11 | 0 |
| 2 | 1 | 7 | 3 | 12 | 0 |
| 3 | 0 | 8 | 0 | 13 | 0 |
| 4 | 3 | 9 | 0 | 14 | 0 |

Now that our groupings are defined, let:

- $\mathbf{Y}_{sj}$ denote the $4 \times 1$ complete data vector for pair $j$ in sequence $s$.
- $\mathbf{Y}_{psj}$ denote the $r \times 1$ observed (reduced) data vector for pair $j$ in pattern $p$ and sequence $s$.
- $\beta^{(g)} = (\mu_{1A}^{(g)}, \mu_{1B}^{(g)}, \mu_{2A}^{(g)}, \mu_{2B}^{(g)}, \rho_1^{(g)}, \rho_2^{(g)}, \nu_1^{(g)}, \nu_2^{(g)})'$ denote the $8 \times 1$ fixed effects location parameter vector for group $g = C, D$, and $P$.
- $\mathbf{X}_{sj}$ denote the $4 \times 8$ design matrix for pair $j$ in sequence $s = 1, 2$ that links the complete data vector $\mathbf{Y}_{sj}$ to $\beta^{(C)}$.
- $\mathbf{\Sigma}$ denote the $4 \times 4$ complete variance-covariance matrix (as defined previously) for the complete data vector $\mathbf{Y}_{sj}$ (notice that we do not assume that the variance-covariance parameters differ across the patterns).
- $\mathbf{E}_{ps}$ be an $r \times 4$ submatrix of the $4 \times 4$ identity matrix for pattern $p$ and sequence $s$, in which rows of the identity matrix are removed according to the missing values in $\mathbf{Y}_{sj}$. The rows of the identity matrix are always considered in 1A, 1B, 2A, 2B order, and $r$ equals the number of non-missing observations.

Then, we define our pattern-mixture model as:

$$\mathbf{Y}_{psj} = \mathbf{E}_{ps}\mathbf{Y}_{sj} \sim N_r\left(\mathbf{E}_{ps}\mathbf{X}_{sj}\beta^{(g)}, \mathbf{E}_{ps}\mathbf{\Sigma}\mathbf{E}'_{ps}\right),$$

where $g = C$ for patterns $p = 0, 10, 11$, and 12; $g = D$ for patterns $p = 1, 2, 3, 6, 7, 13$ and 14; and $g = P$ for patterns $p = 4, 5, 8$ and 9.

For example, pattern $p = 2$ and sequence $s = 1$ is missing the 1B measurement. Therefore, $r = 3$ and:

$$\mathbf{E}_{21} = \begin{pmatrix} 1 & 0 & 0 & 0 \\ 0 & 0 & 1 & 0 \\ 0 & 0 & 0 & 1 \end{pmatrix}.$$

The mean of the response $Y_{21j} = (Y_{21j1A}, Y_{21j2A}, Y_{21j2B})'$ is:

$$\mathbf{E}_{21}\mathbf{X}_{1j}\beta^{(D)} = \begin{pmatrix} \mu_{1A}^{(D)} + \rho_1^{(D)} + \nu_1^{(D)} \\ \mu_{2A}^{(D)} + \rho_2^{(D)} + \nu_2^{(D)} \\ \mu_{2B}^{(D)} - \rho_2^{(D)} + \nu_2^{(D)} \end{pmatrix}$$

and the $3 \times 3$ variance-covariance matrix is:

$$\mathbf{E}_{21}\mathbf{\Sigma}\mathbf{E}'_{21} = \begin{bmatrix} \sigma_{1A,1A} & \sigma_{1A,2A} & \sigma_{1A,2B} \\ \sigma_{1A,2A} & \sigma_{2A,2A} & \sigma_{2A,2B} \\ \sigma_{1A,2B} & \sigma_{2A,2B} & \sigma_{2B,2B} \end{bmatrix}.$$

The proposed means and variances of the remaining patterns $p$ and sequences $s$ can be obtained similarly.

### 2.4. Maximum likelihood estimation

We can readily obtain maximum likelihood (ML) estimates of the parameters of our pattern-mixture model using available software such as PROC MIXED in SAS 9.1. The following is sample code from SAS PROC MIXED that can be used to find the ML estimates of the covariance parameters and the pattern-specific location parameters:

```
PROC MIXED DATA = barge METHOD = ML;
    CLASS pairid position grp;
    MODEL response = mu1A(grp) mu1B(grp) mu2A(grp) mu2B(grp)
                     rho1(grp) rho2(grp) nu1(grp) nu2(grp)/NOINT S;
    REPEATED position / SUBJECT = pairid TYPE = UN;
RUN;
```

The variable *pairid* represents the identification number of the pairs. The variable *position* represents the order of the data in the pair's quadrivariate response vector: 1A is position 1, 1B is position 2, 2A is position 3 and 2B is position 4. The variable *grp* represents the three groups of patterns: for our example, groups C, D, and P.

As an alternative to using SAS PROC MIXED, we can use the log-likelihood function and a matrix programming language, such as S-Plus or SAS/IML, to write our own Fisher scoring or Newton-Raphson algorithm to find the ML estimates. Let:

- $n_{ps}$ denote the number of pairs falling in pattern $p$ of sequence $s$, and $n_s$ denote the number of pairs falling in sequence $s$.
- $\pi_{ps}$ denote the true proportion of pairs falling in pattern $p$ of sequence $s$.
- $\phi_{id}$ denote the set of 62 identifiable parameters – 24 effects $(\mu_{1A}^{(C)}, \ldots, \nu_2^{(P)})$, 28 proportions $(\pi_{01}, \pi_{02}, \ldots, \pi_{13,1}, \pi_{13,2})$, and 10 covariance parameters $(\sigma_{1A,1A}, \ldots, \sigma_{2B,2B})$.

Then, letting $\boldsymbol{\Sigma}_{ps}^* = \mathbf{E}_{ps}\boldsymbol{\Sigma}\mathbf{E}_{ps}'$ and using Lagrange multipliers $\lambda_1$ and $\lambda_2$ with constraints $\pi_{01} + \cdots + \pi_{14,1} = 1$ and $\pi_{02} + \cdots + \pi_{14,2} = 1$, respectively, the log likelihood function is:

$$\log L_{\mathbf{Y},\mathbf{n}}(\phi_{id}) \approx$$
$$\left[\sum_{p=0}^{14}\sum_{s=1}^{2} n_{ps} \log \pi_{ps}\right] - \lambda_1\left(1 - \pi_{01} - \cdots - \pi_{14,1}\right)$$
$$-\lambda_2\left(1 - \pi_{02} - \cdots - \pi_{14,2}\right) \frac{1}{2}\sum_{p=0}^{14}\sum_{s=1}^{2}\sum_{j=1}^{n_{ps}} \log\left|\boldsymbol{\Sigma}_{ps}^*\right|$$
$$-\frac{1}{2}\sum_{p=0}^{14}\sum_{s=1}^{2}\sum_{j=1}^{n_{ps}}\left(\mathbf{E}_{ps}\mathbf{y}_{sj} - \mathbf{E}_{ps}\mathbf{X}_{sj}\beta^{(g)}\right)'\boldsymbol{\Sigma}_{ps}^{*-1}\left(\mathbf{E}_{ps}\mathbf{y}_{sj} - \mathbf{E}_{ps}\mathbf{X}_{sj}\beta^{(g)}\right)$$

where the value of $g$ depends on pattern $p$ ($g = C$ for $p = 0, 10, 11, 12$; $g = D$ for $p = 1, 2, 3, 6, 7, 13, 14$; and $g = P$ for $p = 4, 5, 8, 9$).

### 2.5. *Restricted maximum likelihood estimation*

Just as is the case for ML estimation, it is possible to use SAS PROC MIXED to find restricted maximum likelihood (REML) estimates of the parameters of the pattern mixture model. We only need to make one minor modification to the SAS PROC MIXED code previously used to find the ML estimates (change the "METHOD = ML" option to "METHOD = REML" option).

Again, alternatively we can use the restricted likelihood function and a matrix programming language, such as S-Plus or SAS/IML, to write one's own Fisher scoring or Newton-Raphson algorithm to find the REML estimates.

Let:

- $\mathbf{n} = (n_{01}, \ldots, n_{14,2})'$ be the vector of sample sizes $n_{ps}$ for pattern $p$ and sequence $s$ (with $N = \sum_{p=0}^{14} \sum_{s=1}^{2} n_{ps}$).
- $\pi = (\pi_{01}, \ldots, \pi_{14,2})'$ be the vector of true population proportions for pattern $p$ and sequence $s$.
- $\mathbf{Y} = (\mathbf{y}_{011}, \ldots, \mathbf{y}_{01n_{01}}, \ldots, \mathbf{y}_{14,21}, \ldots, \mathbf{y}_{14,2n_{14,2}})'$ be the $N \times 1$ vector of the entire set of observed responses.
- $\beta = (\beta^{(C)\prime}, \beta^{(D)\prime}, \beta^{(P)\prime})'$ be the $24 \times 1$ vector of fixed effect location parameters.
- $\mathbf{X}$ be the $N \times 24$ full design matrix linking $\mathbf{Y}$ to $\beta$.
- $\mathbf{\Omega}$ be the $N \times N$ block-diagonal variance-covariance matrix of $\mathbf{Y}$.
- $\theta$ be a vector containing the 10 identifiable variance components $(\sigma_{1A,1A}, \ldots, \sigma_{2B,2B})$.

Then the joint likelihood of $\mathbf{Y}$ and $\mathbf{n}$ can be written as the factorization:

$$L_{\mathbf{Y},\mathbf{n}}(\phi_{id}) = L_{\mathbf{Y}|\mathbf{n}}(\beta, \theta) \times L_{\mathbf{n}}(\pi),$$

where $\mathbf{Y}|\mathbf{n}$ is a multivariate $N$-normal with mean $\mathbf{X}\beta$ and variance $\mathbf{\Omega}$, and $\mathbf{n}$ is multinomial with parameters $N$ and $\pi$. Then, the general procedure behind REML estimation of transforming the data vector $\mathbf{Y}|\mathbf{n}$ such that the likelihood $L_{\mathbf{Y}|\mathbf{n}}(\beta, \theta)$ factors into two components can be applied, yielding:

$$\log L_{\mathbf{Y},\mathbf{n}}(\phi_{id}) = \log L'(\theta) + \log L''(\beta, \theta) + \log L_{\mathbf{n}}(\pi).$$

Letting $\mathbf{P} = \mathbf{\Omega}^{-1} - \mathbf{\Omega}^{-1}\mathbf{X}(\mathbf{X}'\mathbf{\Omega}^{-1}\mathbf{X})^{-1}\mathbf{X}'\mathbf{\Omega}^{-1}$, the log REML likelihood function is:

$$\log L'(\theta) \approx -\frac{1}{2}\log|\mathbf{\Omega}| - \frac{1}{2}\log\left|\mathbf{X}'\mathbf{\Omega}^{-1}\mathbf{X}\right| - \frac{1}{2}\mathbf{Y}'\mathbf{P}\mathbf{Y}$$

and:

$$\log L''(\beta, \theta) = -\frac{1}{2}\log\left|(\mathbf{X}'\mathbf{\Omega}^{-1}\mathbf{X})^{-1}\right| - \frac{1}{2}\left(\widehat{\beta} - \beta\right)'(\mathbf{X}'\mathbf{\Omega}^{-1}\mathbf{X})\left(\widehat{\beta} - \beta\right),$$

where $\widehat{\beta}$ has the form of the generalized least squares (GLS) estimator $\widehat{\beta} = (\mathbf{X}'\mathbf{\Omega}^{-1}\mathbf{X})^{-1}\mathbf{X}'\mathbf{\Omega}^{-1}\mathbf{Y}$. Using Lagrange multipliers, $\log L_{\mathbf{n}}(\pi)$ is the log-likelihood of the two independent multinomial samples:

$$\log L_{\mathbf{n}}(\pi) = \left[\sum_{p=0}^{14}\sum_{s=1}^{2} n_{ps} \log \pi_{ps}\right] - \lambda_1 (1 - \pi_{01} - \cdots - \pi_{14,1})$$
$$- \lambda_2 (1 - \pi_{02} - \cdots - \pi_{14,2})$$

## 2.6. Inference about $\gamma$

Thus far, we only have addressed estimation of the pattern-specific parameters, such as $\mu_{1A}^{(C)}, \rho_1^{(D)}$, and $\nu_1^{(P)}$. Interest, however, typically lies in inference about the overall population parameters, such as $\mu_{1A}$, $\rho_1$, and $\nu_1$, not the pattern-specific ones. Therefore, here we extend our work of the previous sections by addressing estimation of the overall population parameters. Furthermore, since use of the paired $2 \times 2$ crossover design will often be motivated by interest in inference about the interaction parameter:

$$\gamma = (\mu_{1A} - \mu_{1B}) - (\mu_{2A} - \mu_{2B}),$$

we also derive the asymptotic distribution of $\widehat{\gamma}$ here. Since the expression for the point estimator $\widehat{\gamma}$ and the form of the asymptotic distribution of $\widehat{\gamma}$ are the same regardless of whether REML or ML estimation is used in obtaining $\widehat{\gamma}$, for the sake of simplicity, we proceed assuming that ML estimation is used.

In order to define point estimators of the overall population parameters, let $\widehat{\pi}_g$ denote the proportion of pairs in group $g$ for $g = C, D$, and $P$. Then, for $k = 1, 2$ and $l = A, B$, the weighted average:

$$
\begin{aligned}
(2.1) \quad \widehat{\mu}_{kl} &= \widehat{\pi}_C \widehat{\mu}_{kl}^{(C)} + \widehat{\pi}_D \widehat{\mu}_{kl}^{(D)} + \widehat{\pi}_P \widehat{\mu}_{kl}^{(P)} \\
&= \widehat{\pi}_C \widehat{\mu}_{kl}^{(C)} + \widehat{\pi}_D \widehat{\mu}_{kl}^{(D)} + (1 - \widehat{\pi}_C - \widehat{\pi}_D) \widehat{\mu}_{kl}^{(P)}
\end{aligned}
$$

is the ML estimator of the overall population treatment-by-genotype mean. And, therefore $\widehat{\gamma} = (\widehat{\mu}_{1A} - \widehat{\mu}_{1B}) - (\widehat{\mu}_{2A} - \widehat{\mu}_{2B})$ is the ML estimator of $\gamma$.

The estimator $\widehat{\gamma}$ is a function of $\widehat{\pi}_g$ and $\widehat{\mu}_{kl}^{(g)}$, the estimated proportions and pattern-specific means. Therefore, we derive the asymptotic distribution of $\widehat{\gamma}$ by first using the asymptotic normality of ML estimators and then by an application of the delta method.

Let:

- $\widehat{\pi} = (\widehat{\pi}_C, \widehat{\pi}_D)'$ be the $2 \times 1$ vector containing the estimated proportion of pairs falling into each of the first two groups.
- $\pi = (\pi_C, \pi_D)'$ be the $2 \times 1$ vector containing the corresponding two true proportions.
- $\mathbf{V}(\widehat{\pi}) = Diag(\pi) - \pi\pi'$ denote the $2 \times 2$ variance matrix of $\widehat{\pi}$.
- $\widehat{\mu} = (\widehat{\mu}_{1A}^{(C)}, \widehat{\mu}_{1A}^{(D)}, \widehat{\mu}_{1A}^{(P)}, \ldots, \widehat{\mu}_{2B}^{(C)}, \widehat{\mu}_{2B}^{(D)}, \widehat{\mu}_{2B}^{(P)})$ be the $12 \times 1$ vector containing the estimated pattern-specific mean parameters.
- $\mu = (\mu_{1A}^{(C)}, \mu_{1A}^{(D)}, \mu_{1A}^{(P)}, \ldots, \mu_{2B}^{(C)}, \mu_{2B}^{(D)}, \mu_{2B}^{(P)})$ be the $12 \times 1$ vector containing the corresponding true pattern-specific mean parameters.
- $\mathbf{V}(\widehat{\mu})$ denote the $12 \times 12$ variance matrix of $\widehat{\mu}$.

Then, by the asymptotic normality of maximum likelihood estimators and by the independence of $\widehat{\pi}$ and $\widehat{\mu}$:

$$
\sqrt{N}\left( \begin{pmatrix} \widehat{\pi}_N \\ \widehat{\mu}_N \end{pmatrix} - \begin{pmatrix} \pi \\ \mu \end{pmatrix} \right)
$$

is asymptotically multivariate normal with mean $\mathbf{0}_{14 \times 1}$ and $14 \times 14$ variance-covariance matrix:

$$
\mathbf{V} = \begin{pmatrix} \mathbf{V}(\widehat{\pi})_{2 \times 2} & \mathbf{0}_{2 \times 12} \\ \mathbf{0}_{12 \times 2} & \mathbf{V}(\widehat{\mu})_{12 \times 12} \end{pmatrix}.
$$

Now, defining the $4 \times 1$ vector of functions $g((\widehat{\pi}', \widehat{\mu}')') = (\widehat{\mu}_{1A}, \widehat{\mu}_{1B}, \widehat{\mu}_{2A}, \widehat{\mu}_{2B})'$ where the functions $\widehat{\mu}_{1A}, \widehat{\mu}_{1B}, \widehat{\mu}_{2A}$ and $\widehat{\mu}_{2B}$ are as defined by (2.1), the $4 \times 14$ Jacobian matrix equals $\mathbf{J} = [\mathbf{J}_1 | \mathbf{J}_2]$ where $\mathbf{J}_1$ is the $4 \times 2$ matrix:

$$
\mathbf{J}_1 = \begin{pmatrix} \widehat{\mu}_{1A}^{(C)} - \widehat{\mu}_{1A}^{(P)} & \widehat{\mu}_{1A}^{(D)} - \widehat{\mu}_{1A}^{(P)} \\ \widehat{\mu}_{1B}^{(C)} - \widehat{\mu}_{1B}^{(P)} & \widehat{\mu}_{1B}^{(D)} - \widehat{\mu}_{1B}^{(P)} \\ \widehat{\mu}_{2A}^{(C)} - \widehat{\mu}_{2A}^{(P)} & \widehat{\mu}_{2A}^{(D)} - \widehat{\mu}_{2A}^{(P)} \\ \widehat{\mu}_{2B}^{(C)} - \widehat{\mu}_{2B}^{(P)} & \widehat{\mu}_{2B}^{(D)} - \widehat{\mu}_{2B}^{(P)} \end{pmatrix}
$$

and $\mathbf{J}_2$ is the $4 \times 12$ matrix:

$$\mathbf{J}_2 = \begin{pmatrix} \widehat{\pi}_C & \widehat{\pi}_D & \widehat{\pi}_P & 0 & 0 & 0 & 0 & 0 & 0 & 0 & 0 & 0 \\ 0 & 0 & 0 & \widehat{\pi}_C & \widehat{\pi}_D & \widehat{\pi}_P & 0 & 0 & 0 & 0 & 0 & 0 \\ 0 & 0 & 0 & 0 & 0 & 0 & \widehat{\pi}_C & \widehat{\pi}_D & \widehat{\pi}_P & 0 & 0 & 0 \\ 0 & 0 & 0 & 0 & 0 & 0 & 0 & 0 & 0 & \widehat{\pi}_C & \widehat{\pi}_D & \widehat{\pi}_P \end{pmatrix}.$$

Then, an application of the delta method yields that:

$$\sqrt{N} \left( \begin{pmatrix} \widehat{\mu}_{1A} \\ \widehat{\mu}_{1B} \\ \widehat{\mu}_{2A} \\ \widehat{\mu}_{2B} \end{pmatrix} - \begin{pmatrix} \mu_{1A} \\ \mu_{1B} \\ \mu_{2A} \\ \mu_{2B} \end{pmatrix} \right)$$

is asymptotically multivariate normal with mean $\mathbf{0}_{4 \times 1}$ and $4 \times 4$ variance matrix $\mathbf{JVJ}'$.

Now, we can construct any linear combination of the population mean parameters as $\gamma = \mathbf{c}'(\mu_{1A}, \mu_{1B}, \mu_{2A}, \mu_{2B})'$, but the contrast vector of particular interest in our application is:

$$\mathbf{c} = \begin{pmatrix} 1 & -1 & -1 & 1 \end{pmatrix}'.$$

Then

$$\widehat{\gamma} = g\left( \begin{pmatrix} \widehat{\mu}_{1A} & \widehat{\mu}_{1B} & \widehat{\mu}_{2A} & \widehat{\mu}_{2B} \end{pmatrix}' \right) = \mathbf{c}'(\widehat{\mu}_{1A}, \widehat{\mu}_{1B}, \widehat{\mu}_{2A}, \widehat{\mu}_{2B})'$$

and

$$\frac{\partial \widehat{\gamma}}{\partial \begin{pmatrix} \widehat{\mu}_{1A} & \widehat{\mu}_{1B} & \widehat{\mu}_{2A} & \widehat{\mu}_{2B} \end{pmatrix}} = \mathbf{c}'.$$

Therefore, one final application of the delta method yields that $\sqrt{N}(\widehat{\gamma} - \gamma)$ is asymptotically normal with mean 0 and variance $V(\widehat{\gamma}) = \mathbf{c}'\mathbf{JVJ}'\mathbf{c}$. The resulting asymptotic variance $V(\widehat{\gamma})$ depends on $\widehat{\pi}_g$, $\widehat{\mu}_{kl}^{(g)}$, $\mathbf{V}(\widehat{\pi})$ and $\mathbf{V}(\widehat{\mu})$, and therefore must be estimated.

## 3. Results

As described in the Introduction, the BARGE Study was a randomized paired $2 \times 2$ crossover trial in asthma patients comparing the effects of regularly-scheduled inhaled albuterol to placebo (Israel et al. [6]). Patients who met the asthma eligibility criteria were genotyped at the 16th position of the beta-agonist receptor gene (Arg/Arg, Gly/Gly, Arg/Gly). The heterozygotes (Arg/Gly) were excluded from the study. The Arg/Arg and Gly/Gly patients, on the other hand, were matched according to their pulmonary function (as determined by their forced expiratory volume in one second) and randomized together to the same treatment sequence. The primary outcome variable was the difference in the patients' morning peak expiratory flow rates (AM PEFR, measured in liters per minute) at the start and end of the 16-week treatment periods. A change of 25 liters per minute was the effect size considered to be clinically meaningful for the sample size calculation made for the BARGE study protocol.

The BARGE Study consisted of 78 randomized patients, of whom 71 yielded data for analysis – 35 with the Arg/Arg genotype and 36 with the Gly/Gly genotype at the 16th position of the beta-agonist receptor gene. There were 40 sets of subjects, and as mentioned in Section 2.3, $n_C = 29$, $n_D = 6$, and $n_P = 5$. Because of the small samples sizes for groups D and P, which yielded non-estimable standard errors for some of their estimated parameters, these two groups were pooled. Thus, the final pattern-mixture model for the analysis of the BARGE Study data consisted of only two groups (C and D+P).

A general linear model with correlated errors, as described in Sections 2.4 and 2.5, was applied to group C and group D+P separately that included effects for period, sequence, and treatment. Because group D+P had a small sample size, a common $4 \times 4$ variance-covariance matrix was assumed for group C and group D+P. The REML estimates for the treatment effects, along with their standard errors, appear in Table 5. It is not surprising that the standard errors for the treatment effect estimates in group D+P are much larger than those of group C because (1) group C consisted of 29 sets of patients with 116 observations, and (2) group D+P consisted of 11 sets of patients with 20 observations. The REML estimates for the weighted averages of the treatment effects, along with standard errors, appear in Table 6, based on the pattern-mixture methodology described in the previous section.

The estimate of $\gamma = (\mu_{RA} - \mu_{RP}) - (\mu_{GA} - \mu_{GP})$, the genotype × treatment interaction term, is –11.7 liters per minute with standard error 22.3 ($p = 0.60$). Thus, the analysis based on the proposed pattern-mixture model does not yield a statistically significant result. If the patterns (groups C and D+P) are ignored in the analysis, then the estimate of $\gamma$ is $-15.8$ liters per minute with standard error 19.4 ($p = 0.42$). Although the results of the two analyses (with and without pattern-mixture modeling) yield similar scientific conclusions, the analysis based on the proposed pattern-mixture model is slightly more conservative.

In the BARGE scientific manuscript, Israel et al. [6] applied a mixed-effects linear model to the repeated measurements data within each of the 16-week treatment periods. In other words, they used all of the available data, not just the patients' final 16-week measurement minus their baseline measurement. Such an approach yielded a more powerful and sensitive analysis than the analysis described here because many of the patients in pattern-mixture group D had partial data during the 16-week treatment periods which contributed to the overall analysis. Indeed, the authors of the BARGE manuscript estimated the genotype × treatment interaction term as 24.0 liters per minute with standard error 6.3 ($p < 0.01$). In a future communication, we plan to develop the pattern-mixture version of the mixed-effects linear model for the paired $2 \times 2$ crossover design that will accommodate repeated measurements within the randomized treatment periods.

TABLE 5
*The REML estimates of the treatment effects from the mixed-effects model analysis of groups C and D+P separately. Arg/Arg = R, Gly/Gly = G, Albuterol = A, Placebo = P*

| Group | Genotype | Treatment | Change in AM PEFR $\widehat{\mu}_{kl}^{(g)} =$ Mean (Std Error) |
|---|---|---|---|
| C | Arg/Arg | Placebo | $\widehat{\mu}_{RP}^{(C)} = 20.4\ (11.7)$ |
| C | Arg/Arg | Albuterol | $\widehat{\mu}_{RA}^{(C)} = 8.1\ (10.3)$ |
| C | Gly/Gly | Placebo | $\widehat{\mu}_{GP}^{(C)} = 12.6\ (17.7)$ |
| C | Gly/Gly | Albuterol | $\widehat{\mu}_{GA}^{(C)} = 22.3\ (20.3)$ |
| D+P | Arg/Arg | Placebo | $\widehat{\mu}_{RP}^{(D+P)} = -23.7\ (35.7)$ |
| D+P | Arg/Arg | Albuterol | $\widehat{\mu}_{RA}^{(D+P)} = 12.0\ (40.1)$ |
| D+P | Gly/Gly | Placebo | $\widehat{\mu}_{GP}^{(D+P)} = -66.8\ (45.5)$ |
| D+P | Gly/Gly | Albuterol | $\widehat{\mu}_{GP}^{(D+P)} = -46.4\ (45.4)$ |

TABLE 6
*The REML estimates of the weighted averages of the treatment effects from the mixed-effects model analysis. Arg/Arg = R, Gly/Gly = G, Albuterol = A, Placebo = P*

| Genotype | Treatment | Change in AM PEFR Mean (Std Error) |
|---|---|---|
| Arg/Arg | Placebo | $\widehat{\mu}_{RP} = 8.3\ (13.0)$ |
| Arg/Arg | Albuterol | $\widehat{\mu}_{RA} = 9.2\ (13.3)$ |
| Gly/Gly | Placebo | $\widehat{\mu}_{GP} = -9.2\ (17.9)$ |
| Gly/Gly | Albuterol | $\widehat{\mu}_{GA} = 3.4\ (19.3)$ |

## 4. Discussion

We have demonstrated the development of a pattern-mixture version of a general linear model with correlated errors for the paired $2 \times 2$ crossover design. As just mentioned, in a future communication we will extend our work described herein to account for repeated measurements within the randomized treatment periods. Such an approach will optimize the pattern mixture approach because dropouts could be redefined as early dropouts or late dropouts.

A weakness of the pattern mixture approach to the analysis of data from a paired $2 \times 2$ crossover design is that one or more of the pattern groups could represent a reasonable proportion of the sample, yet provide little data for analysis. Case in point, this occurred with the BARGE Study data presented here. Eleven of the 40 sets (27.5%) were categorized into group D+P, yet group D+P provided only 20 of the 136 observations (14.7%) in the data set. Thus, the influence of group D+P on the results of the analysis is disproportionate to the actual number of observations that it contributes to the analysis. This weakness could be minimized by using available repeated measurements within the randomized treatment periods.

One popular alternative to the pattern-mixture model approach is multiple data imputation (e.g., see Schafer [16]). Multiple data imputation involves the construction of a complete data set via estimation of the missing values according to an appropriate probability model, repeating this $m$ times, estimating the treatment effects within each data set, and averaging the effects across the $m$ data sets. Multiple data imputation affords two advantages: (1) it provides complete data sets that are more easily analyzed via standard statistical methods, and (2) it appropriately accounts for the variability in estimating the missing values. The disadvantage of multiple data imputation is that the probability model for generating imputed values could be misspecified. If so, then the combined results could be biased. For example, the probability model may be based on the estimated mean and variance structure of those experimental units with complete observations. If the data from dropouts have an inherently different mean and variance structure, then the imputed values could be misrepresentative.

**Acknowledgments.** The authors wish to thank an anonymous referee, whose valuable comments improved the presentation of this material.